\documentclass{amsproc}
\usepackage[mathscr]{eucal}
\usepackage{amsthm}
\usepackage{amscd}
\usepackage{pb-diagram}

\def\Bbb{\mathbb}
\def\frak{\mathfrak}

\newenvironment{pf*}[1]{\proof[#1]}{\endproof}
\newcommand{\rom}{\textup}
\newtheorem{Theorem}[equation]{Theorem}

\newtheorem{Lemma}[equation]{Lemma}
\newtheorem{Proposition}[equation]{Proposition}

\newtheorem{Conjecture}[equation]{Conjecture}

\theoremstyle{definition}
\newtheorem{Definition}[equation]{Definition}
\newtheorem{Example}[equation]{Example}

\theoremstyle{remark}
\newtheorem{Remark}[equation]{Remark}

\numberwithin{equation}{subsection}
\newcommand{\thmref}[1]{Theorem~\ref{#1}}
\newcommand{\secref}[1]{\S\ref{#1}}
\newcommand{\lemref}[1]{Lemma~\ref{#1}}

\newcommand{\subsecref}[1]{\S\ref{#1}}

\newcommand{\defeq}{\overset{\operatorname{\scriptstyle def.}}{=}}
\newcommand{\C}{{\Bbb C}}
\newcommand{\Z}{{\Bbb Z}}

\newcommand{\R}{{\Bbb R}}
\newcommand{\GL}{\operatorname{GL}}

\newcommand{\algsl}{\operatorname{\frak{sl}}} % because \sl="slant"

\newcommand{\g}{{\frak g}}

\newcommand{\Ker}{\operatorname{Ker}}

\newcommand{\id}{\operatorname{id}}
\newcommand{\ve}{\varepsilon}
\newcommand{\M}{{\frak M}} % the moduli space
\newcommand{\Mreg}{{\frak M_0^{\operatorname{reg}}}}
\newcommand{\La}{{\frak L}} % the Lagrangian variety
\newcommand{\bv}{{\mathbf v}} % the vector v
\newcommand{\bw}{{\mathbf w}} % the vector w
\newcommand{\topdeg}{\operatorname{top}} % top degree
\newcommand{\Uq}{{\mathbf U}_q(\mathfrak g)} % the QUE algebra

\newcommand{\Ua}{{\mathbf U}_q(\widehat{\mathfrak g})} % the quantum affine
\newcommand{\Ul}{{\mathbf U}_q({\mathbf L}{\mathfrak g})} % the quantum
\newcommand{\Ue}{{\mathbf U}_{\varepsilon}(\mathfrak g)}

\newcommand{\Uli}{{\mathbf U}^{\Z}_q({\mathbf L}\mathfrak g)}
\newcommand{\Ule}{{\mathbf U}_{\varepsilon}({\mathbf L}\mathfrak g)}
\newcommand{\bU}{\mathbf U} % the enveloping algebra
\newcommand{\Wedge}{{\textstyle \bigwedge}}
\newcommand{\Zw}{Z(\bw)}
\newcommand{\Mw}{\M(\bw)}
\newcommand{\Law}{\La(\bw)}

\newcommand{\Lg}{\mathbf L\g}

\begin{document}
\title[$t$-analogue of the $q$-characters]
{$t$-analogue of the $q$-characters of finite dimensional
representations of quantum affine algebras
%\\ {\rm \small (Preliminary version: \today)}
}
\author{Hiraku Nakajima}
\address{Department of Mathematics, Kyoto University, Kyoto 606-8502,
Japan
%\\Email:nakajima@kusm.kyoto-u.ac.jp}
}
\email{nakajima@kusm.kyoto-u.ac.jp}
\urladdr{http://www.kusm.kyoto-u.ac.jp/\textasciitilde nakajima}
\thanks{Supported by the Grant-in-aid
for Scientific Research (No.11740011), the Ministry of Education,
Japan.}
%
%\subjclass{Primary 17B37;
%Secondary 14D21, 14L30, 16G20, 33D80}
%
\begin{abstract}
%\abstracts{
Frenkel-Reshetikhin introduced $q$-characters of finite dimensional
representations of quantum affine algebras~\cite{FR}.
We give a combinatorial algorithm to compute them for all simple
modules.
Our tool is $t$-analogue of the $q$-characters, which is similar to
Kazhdan-Lusztig polynomials, and our algorithm has a resemblance with
their definition.

We need the theory of quiver varieties for the definition of
$t$-analogues and the proof. But it appear only in the last section.
The rest of the paper is devoted to an explanation of the algorithm,
which one can read without the knowledge about quiver varieties.
A proof is given only in part. A full proof will appear elsewhere.
\end{abstract}
%}
\maketitle

\section{The quantum loop algebra}

Let $\g$ be a simple Lie algebra of type ADE over $\C$, $\Lg = \g
\otimes \C[z,z^{-1}]$ be its loop algebra, and $\Ul$ be its quantum
universal enveloping algebra, or the quantum loop algebra for short.
It is a subquotient of the quantum affine algebra $\Ua$, i.e., without
central extension and degree operator.
Let $I$ be the set of simple roots, $P$ be the weight lattice, and
$P^*$ be its dual lattice (all for $\g$).
The algebra has the so-called Drinfeld's new realization: It is a
$\C(q)$-algebra with generators $q^h$, $e_{k,r}$, $f_{k,r}$, $h_{k,n}$
($h\in P^*$, $k\in I$, $r\in\Z$, $n\in \Z\setminus \{0\}$) with
certain relations (see e.g., \cite[12.2]{CP-book}).

The algebra $\Ul$ is a Hopf algebra, where the coproduct is defined
using the Drinfeld-Jimbo realization of $\Ul$. So a tensor product
$M\otimes_{\C(q)} M'$ of $\Ul$-modules $M$, $M'$ has a structure of a
$\Ul$-module.

Let $\Ule$ be its specialization at $q = \ve\in\C^*$.
For precise definition of the specialization, we first introduce an
integral form $\Uli$ of $\Ul$ and set $\Ule =
\Uli\otimes_{\Z[q,q^{-1}]} \C$, where $\Z[q,q^{-1}]\to \C$ is given by
$q^{\pm 1}\mapsto \ve^{\pm 1}$. See \cite{CP-roots} for detail. But we
assume $\ve$ is {\it not\/} a root of unity in this paper. So we just
replace $q$ by $\ve$ in the definition of $\Ul$.

The quantum loop algebra $\Ul$ contains the quantum enveloping algebra 
$\Uq$ for the finite dimensional Lie algebra $\g$ as a subalgebra. The 
specialization $\Ule$ contains the specialization $\Ue$ of $\Uq$.

\subsection{Finite dimensional representations of $\Ule$}

The algebra $\Ule$ contains a commutative subalgebra generated by
$q^h$, $h_{k,n}$ ($h\in P^*$, $k\in I$, $n\in\Z\setminus \{0\}$).
Let us introduce generating functions $\psi_k^\pm(z)$ ($k\in I$) by
\begin{equation*}
   \psi^{\pm}_k(z)
%  = \sum_{r=0}^\infty \psi^{\pm}_{k,\pm r} z^{\mp r}
  \defeq q^{\pm h_k}
   \exp\left(\pm (q-q^{-1})\sum_{m=1}^\infty h_{k,\pm m} z^{\mp m}\right).
\end{equation*}

A $\Ule$-module $M$ is called of {\it type $1$\/} if $M$ has a weight space
decomposition as a $\Ue$-module:
\begin{equation*}
  M = \bigoplus_{\lambda\in P} M(\lambda), \qquad
  M(\lambda) = \left\{ m\in M \left|\, q^h\ast m 
  = \ve^{\langle h,\lambda\rangle} m\right.\right\}.
\end{equation*}
We will only consider type $1$ modules in this paper.

A type $1$ module $M$ is an {\it l-highest weight module\/} ('{\it
l\/}' stands for the loop) if there exists a vector $m_0\in M$ such
that
\begin{gather*}
  e_{k,r}\ast m_0 = 0, \qquad \Ule^- \ast m_0 = M, \\
  \psi^{\pm}_{k}(z)\ast m_0 = \Psi^{\pm}_{k}(z) m_0
    \quad\text{for $k\in I$}
\end{gather*}
for some $\Psi^{\pm}_k(z)\in \C[[z^\mp]]$. The pair of
the $I$-tuple
$(\Psi^+(z),\Psi^-(z)) = (\Psi^+_k(z),\linebreak[0]
 \Psi^-_k(z))_{k\in I}\in
(\C[[z^\mp]]^I)^2$ is called the {\it l-highest weight\/} of $M$, and
$m_0$ is called the {\it l-highest weight vector}.

\begin{Theorem}[Chari-Pressley~\cite{CP-rep}]
\rom{(1)} Every finite-dimensional simple \linebreak[4]
$\Ule$-module of type $1$ is 
an {\it l\/}-highest weight module, and its {\it l\/}-highest weight
is given by
\begin{equation}\label{eq:hwt}
   \Psi^\pm_k(z) = \varepsilon^{\deg P_k}
  \left(\frac{P_k(\varepsilon^{-1}/z)}{P_k(\varepsilon/z)}\right)^{\pm}
\end{equation}
for some polynomials $P_k(u)\in\C[u]$ with $P_k(0) = 1$.
Here $\left(\ \right)^{\pm}
\in\C[[z^{\mp}]]$ denotes the expansion at $z = \infty$ and $0$
respectively.

\rom{(2)} Conversely, for given $P_k(u)$ as above, there exists a
finite-dimensional simple {\it l\/}-highest weight $\Ule$-module $M$
of type $1$ such that the {\it l\/}-highest weight is given by the
above formula.

Assigning to $M$ the $I$-tuple $P = (P_k)_{k\in I}\in \C[u]^I$
\rom($P_k(0) = 1$\rom) defines a bijection between the set of all
$P$'s and the set of isomorphism classes of finite-dimensional simple
$\Ule$-modules of type $1$.
\end{Theorem}

We denote by $L_P$ the simple $\Ule$-module associated to $P$.
We call $P$ the {\it Drinfeld polynomial}.
For the abuse of terminology, we also say `$P$ is the {\it
l\/}-highest weight of $L_P$'.

Since $\C\langle q^h$, $h_{k,n}\rangle$ is a commutative subalgebra
of $\Ule$, any $\Ule$-module $M$ decomposes into a direct sum $M =
\bigoplus M(\Psi^+,\Psi^-)$ of generalized eigenspaces, where
\begin{equation*}
\begin{split}
   & M(\Psi^+,\Psi^-) \\
   \defeq\; &  \left\{ m\in M \left|\,
     \text{$(\psi_k^\pm(z) - \Psi_k^\pm(z)\operatorname{Id})^N \ast m
     = 0$ for $k\in I$ and sufficiently large $N$}
   \right\}\right. ,
\end{split}
\end{equation*}
for $\Psi_k^\pm(z)\in\C[[z^\mp]]$. The pair of the $I$-tuple
$(\Psi^+,\Psi^-) = (\Psi^+_k, \Psi^-_k)_{k\in I}$ is called an
{\it l-weight}, and $M(\Psi^+,\Psi^-)$ is called an {\it l-weight
space\/} of $M$ if $M(\Psi^+,\Psi^-) \neq 0$.

\begin{Theorem}[Frenkel-Reshetikhin~\cite{FR}]
Any {\it l\/}-weight of any finite dimensional $\Ule$-module $M$ of
type $1$ has the following form
\begin{equation}\label{eq:lwt}
   \Psi_k^\pm(z) = \varepsilon^{\deg Q_k - \deg R_k}
   \left(\frac{Q_k(\varepsilon^{-1}/z) R_k(\varepsilon/z)}
     {Q_k(\varepsilon/z)R_k(\varepsilon^{-1}/z)}\right)^{\pm}
\end{equation}
for some polynomials
\[
  Q_k(u) = \prod_{i=1}^{s_k} (1 - a_{ki} u), \quad
  R_k(u) = \prod_{j=1}^{r_k} (1-b_{kj} u). 
\]
\end{Theorem}
Again for the abuse of terminology, we also say `$Q/R$ is an {\it
l\/}-weight of $M$'.
We denote the {\it l\/}-weight space $M(\Psi^+,\Psi^-)$ by $M(Q/R)$.

Frenkel-Reshetikhin~\cite{FR} defined the $q$-character of $M$ by
\begin{equation*}
  \chi_q(M) \defeq
  \sum_{Q/R}
  \dim M(Q/R)\; \prod_{k\in I}
  \prod_{i=1}^{s_k}\prod_{j=1}^{r_k} Y_{k,a_{ki}} Y_{k,b_{kj}}^{-1}.
\end{equation*}

\begin{Theorem}[Frenkel-Reshetikhin~\cite{FR}]\label{thm:FR}
\rom{(1)} $\chi_q$ defines an injective ring homomorphism from the
Grothendieck ring $\operatorname{Rep}\Ule$ of finite dimensional
$\Ule$-modules of type $1$ to $\Z[Y_{k,a}^\pm]_{k\in I,a\in\C^*}$ \rom(a
ring of Laurent polynomials in infinitely many variables\rom).

\rom{(2)} If we compose a map $Y_{k,a}^\pm \mapsto y_k^\pm$ \rom(forgetting
`spectral parameters'\rom), it gives the usual character of the
restriction of $M$ to a $\Ue$-module.
\end{Theorem}

\begin{Definition}
A monomial $\displaystyle \prod_{k\in I}
\prod_{i=1}^{s_k}\prod_{j=1}^{r_k} Y_{k,a_{ki}} Y_{k,b_{kj}}^{-1}$
appearing in the $q$-char\-ac\-ter $\chi_q$ is called {\it l-dominant\/}
if $r_k = 0$ for all $k$, i.e., a product of positive powers
of $Y_{k,c}$'s or $1$.
\end{Definition}

If $L_P$ is the simple $\Ule$-module with {\it l\/}-highest weight
$P$, its $q$-character contains an {\it l\/}-dominant monomial
corresponding to the {\it l\/}-highest weight. We denote it by
$m_P$. Its coefficient in $\chi_q(L_P)$ is $1$.

Since $\{ L_P \}_P$ forms a basis of $\operatorname{Rep}\Ule$, we have 
the following useful condition for the simplicity of a finite
dimensional $\Ule$-module $M$ of type $1$:
\begin{equation}
\label{eq:suff}
\text{
If $\chi_q(M)$ contains only one {\it l\/}-dominant term, then
$M$ is simple.
}
\end{equation}

\subsection{Example}
We give examples of $q$-characters.

If $\g = A_n$, we have an evaluation homomorphism
$\operatorname{ev}_a\colon \Ule \to \Ue$ corresponding to $\Lg \to
\g$; $z\mapsto a$ (Jimbo). Hence pullbacks of simple $\Ue$-modules are 
simple $\Ule$-modules.

\begin{Example}
Let $\g = A_1 = \algsl_2$ and $V$ be the $2$-dimensional simple
$\Ule$-module. Then the $q$-character of $M_a =
\operatorname{ev}_a(M)$ is given by\footnote{This can be checked
directly. But it also follows from \thmref{thm:ind} below.}
\begin{equation*}
   \chi_q(M_a)
   = Y_{1,a} + Y_{1,a\ve^2}^{-1}.
\end{equation*}
Since $\chi_q$ is a ring homomorphism, we have
\begin{equation*}
\begin{split}
   \chi_q(M_a\otimes M_b) &=
   \left(Y_{1,a} + Y_{1,a\ve^2}^{-1}\right)
   \left(Y_{1,b} + Y_{1,b\ve^2}^{-1}\right)\\
   &=
   Y_{1,a}Y_{1,b} + Y_{1,a\ve^2}^{-1}Y_{1,b} + Y_{1,a}Y_{1,b\ve^2}^{-1}
   + Y_{1,a\ve^2}^{-1}Y_{1,b\ve^2}^{-1}.
\end{split}
\end{equation*}
If $b\neq a\ve^2, a\ve^{-2}$, then $M_a\otimes M_b$ is simple by the
criterion~\eqref{eq:suff}.

If $b = a\ve^2$ or $a\ve^{-2}$, then the second or third term becomes
$1$. In fact, it is known that $M_a\otimes M_{a\ve^2}$ decomposes (in
$\operatorname{Rep}\Ule$) to a sum $M'_a\oplus M''$, where $M'$ is the
$3$-dimensional simple $\Ule$-module, and $M''$ is the trivial
module. Thus we have
\begin{equation*}
\begin{split}
  \chi_q(M_a\otimes M_{a\ve^2}) 
  & = \chi_q(M'_a) + \chi_q(M'') \\
  & = Y_{1,a}Y_{1,a\ve^2} + Y_{1,a}Y_{1,a\ve^4}^{-1}
   + Y_{1,a\ve^2}^{-1}Y_{1,a\ve^4}^{-1} + 1.
\end{split}
\end{equation*}
See also Examples~\ref{exm:A1''}, \ref{exm:A1}, \ref{exm:A1'}.
\end{Example}

\section{Standard modules}

\subsection{}

In \cite{Na-qaff} we defined a family of finite dimensional
$\Ule$-modules of type $1$ and called them {\it standard
modules}. They are parametrized by the $I$-tuples $P = (P_k)_{k\in
I}\in \C[u]^I$ exactly as simple modules. We denote by $M_P$
associated to $P$. The definition will be recalled in
\secref{sec:quiver}, but we give here their algebraic
identification due to Varagnolo-Vasserot~\cite{VV-std}.

\begin{Definition}
We say $L_P$ an {\it l-fundamental representation\/} if
\begin{equation*}
  P_k(u) =
  \begin{cases}
    1 - su  & \text{if $k = k_0$,} \\
    1 & \text{otherwise},
  \end{cases}
\end{equation*}
for some $s\in\C^*$ and $k_0\in I$.
We denote $L_P$ by $L(\Lambda_{k_0})_s$.
($\Lambda_k$ is the $k$-th fundamental weight of $\g$.)
\end{Definition}

For $s\in\C^*$ and a finite sequence
$(k_\alpha)_\alpha = (k_1, k_2, \dots)$ in $I$ and
a sequence $(n_\alpha)_\alpha = (n_1\ge n_2\ge \dots)$
of integers, we set
\begin{equation*}
   M(s; (k_\alpha)_\alpha, (n_\alpha)_\alpha) \defeq
   L(\Lambda_{k_1})_{\ve^{n_1}s}\otimes
   L(\Lambda_{k_2})_{\ve^{n_2}s}\otimes\cdots.
\end{equation*}
Note that $\Ule$ is {\it not\/} cocommutative Hopf algebra, so the
tensor product depends on the ordering of factors.

\begin{Theorem}[Varagnolo-Vasserot~\cite{VV-std}]\label{thm:VV}
\rom{(1)} A standard module $M$ is isomorphic to a module of the form
\begin{equation*}
\begin{split}
   & \bigotimes_{i} 
   M(s^i; (k^i_{\alpha_i})_{\alpha_i}, (n^i_{\alpha_i})_{\alpha_i}) \\
   =\; &
   M(s^1; (k^1_{\alpha_1})_{\alpha_1}, (n^1_{\alpha_1})_{\alpha_1})
   \otimes
   M(s^2; (k^2_{\alpha_2})_{\alpha_2}, (n^2_{\alpha_2})_{\alpha_2})
   \otimes
   \cdots
   \quad \text{\rom(finite tensor product\rom)}
\end{split}
\end{equation*}
such that $s^i/s^j\notin \ve^{\Z}$ for $i\neq j$ and
$n^i_1 \ge n^i_2 \ge \dots$ for each $i$.

\rom{(2)} The above tensor product is independent of the ordering of
the factors
\linebreak[4]
$M(s^i; (k^i_{\alpha_i})_{\alpha_i},(n^i_{\alpha_i})_{\alpha_i})$.

\rom{(3)} The $I$-tuple of polynomials $P$ corresponding to $M$ is the 
product of Drinfeld polynomials of {\it l\/}-fundamental representations
appearing as factors of $M$.
\end{Theorem}

Note that if $P$ is given, we can define a module $M$ of the above
form by decomposing $P$ into a product of Drinfeld polynomials of {\it
l\/}-fundamental representations. Thus we may denote the above module
by $M_P$.

The following properties of $M_P$ were shown in \cite{Na-qaff}:
\begin{enumerate}
\item $\{ M_P\}$ is a basis of $\operatorname{Rep}\Ule$.
\item $M_P$ is an {\it l\/}-highest weight module with {\it
l\/}-highest weight $P$ (i.e., given by \eqref{eq:hwt}).
\item $L_P$ is the unique simple quotient of $M_P$.
\item $M_P$ depends `continuously' on $P$ in a certain sense.
For example, $\dim M_P$ is independent of $P$.
\item for a generic $P$, $M_P \cong L_P$.
\end{enumerate}

Conjecturally $M_P$ is isomorphic to the specialization of the module
$V^{\max}(\lambda)$, introduced by Kashiwara~\cite{Kas}, and further
studied by Chari-Pressley~\cite{CP:Weyl}.

\section{$t$-analogues of $q$-characters}

A main tool in this paper is a $t$-analogue of the $q$-character:
\begin{equation*}
   \chi_{q,t}\colon \operatorname{Rep}\Ule \to
   \Z[t,t^{-1}][Y_{k,a}^\pm]_{k\in I,a\in\C^*}.
\end{equation*}
This is a homomorphism of additive groups, not of rings, and has the
property $\chi_{q,t=1} = \chi_q$.
We define $\chi_{q,t}$ for all standard modules $M_P$. Since 
$\{ M_P\}_P$ is a basis of $\operatorname{Rep}\Ule$, we can extend it 
linearly to any finite dimensional $\Ule$-modules.

For the definition we need geometric constructions of standard
modules, so we will postpone it to \subsecref{subsec:qchar}. We give
an alternative definition, which is conjecturally the same as the
geometric definition.

\subsection{A conjectural definition}

Let $M = M_P$ be a standard module, $Q/R$ be an {\it l\/}-weight of
$M$, $M(Q/R)$ be the corresponding {\it l\/}-weight space.
Define a filtration on $M(Q/R)$ by
\begin{gather*}
   0 = M^{-1}(Q/R) \subset M^0(Q/R)
   \subset M^1(Q/R) \subset \cdots \\
   M^n(Q/R)
   \defeq \bigcap_k \Ker(\psi_k^\pm(z)-\Psi_k^\pm(z)\id)^{n+1}.
\end{gather*}

\begin{Conjecture}
The $t$-analogue $\chi_{q,t}(M_P)$, defined geometrically in
\subsecref{subsec:qchar}, is equal to
\begin{equation*}
  \chi_{q,t}(M_P)
  = \sum_{Q/R}
  \sum_n t^{2n-d(Q/R,P)}
  \dim \left(M^{n}(Q/R)/M^{n-1}(Q/R)\right) \; m_{Q/R},
\end{equation*}
where $d(Q/R,P)$ is an integer \rom(determined explicitly from $Q/R$,
$P$ by \eqref{eq:dim} below\rom), and $m_{Q/R}$ is a monomial in
$Y_{k,a}^\pm$ corresponding to the {\it l\/}-weight space $M(Q/R)$.
\end{Conjecture}

This definition makes sense for any finite dimensional modules, but is
{\it not\/} well-defined on the Grothendiek group
$\operatorname{Rep}\Ule$. Thus the above does not hold for simple modules.

\subsection{}
A main result of this paper is a combinatorial algorithm for computing
$\chi_{q,t}(M_P)$ and $[M_P:L_Q]$. It is divided into three steps:
\begin{description}
\item[Step 1] Compute $\chi_{q,t}$ for all {\it l\/}-fundamental
repsentations.
\item[Step 2] Compute $\chi_{q,t}(M_P)$ for all standard modules
$M_P$.
\item[Step 3] Express the multiplicity $[M_P: L_Q]$ in terms of
$\chi_{q,t}(M_R)$ for various $R$.
\end{description}

Step 1 is a modification of Frenkel-Mukhin's algorithm~\cite{FM}
for computing $\chi_q$ of {\it l\/}-fundamental representations.
Step 2 is nothing but a study of $\chi_{q,t}$ of tensor products of
{\it l\/}-fundamental representations. Although $\chi_{q,t}$ is not a
ring homomorphism, $\chi_{q,t}$ of tensor products is given by a
simply modified multiplication.
For the proof we use an idea in \cite{Na-hom}.
Step 3 was essentially done in \cite{Na-qaff}.

\section{Step 3}

We start with Step 3. The algorithm is similar to the definition of
Kazhdan-Lusztig polynomials \cite{KL}. It is also similar to the
algorithm for computing the transition matrix between the
canonical basis and the PBW basis of type ADE \cite{Lu:can}.

\subsection{}
Let 
\[
   A_{k,a} \defeq Y_{k,a\ve} Y_{k,a\ve^{-1}}
   \prod_{l:l\neq k} Y_{l,a}^{c_{kl}},
\]
where $c_{kl}$ is the $(k,l)$-entry of the Cartan matrix.

\begin{Definition}
(1) Let $m$, $m'$ be monomials in $Y_{k,a}^\pm$ ($k\in I$,
$a\in\C^*$). We define an ordering $\le$ among monomials by
\begin{equation*}
   m\le m' \Longleftrightarrow
   \text{$\frac{m'}{m}$ is a monimial in $A_{k,a}^{-1}$ ($k\in I$,
   $a\in\C^*$)}.
\end{equation*}
Here a monomial in $A_{k,a}^{-1}$ means a product of nonnegative
powers of $A_{k,a}^{-1}$. It does not contain any factors $A_{k,a}$.

(2) If $\Psi^\pm, \Psi^{\prime\pm}$ are {\it l\/}-weights of finite
dimensional $\Ule$-modules, or $Q/R, Q'/R'$ are related to {\it
l\/}-weights by \eqref{eq:lwt}, we write
$\Psi^\pm\le \Psi^{\prime\pm}$, $Q/R \le Q'/R'$ if the corresponding
monomials $m$, $m'$ satisfy $m\le m'$.
\end{Definition}

Recall that $\chi_q(L_P)$ contains an {\it l\/}-dominant monomial
$m_P$ corresponding to the highest weight vector.
It is known that any monomial $m$ appearing $\chi_q(L_P)$,
$\chi_q(M_P)$ satisfies $m\le m_P$ (\cite[4.1]{FM},
\cite[13.5.2]{Na-qaff}).

Let
\begin{equation*}
   c_{QP}(t) \defeq \text{the coefficient of $m_Q$ in $\chi_{q,t}(M_P)$}.
\end{equation*}
Then $(c_{QP}(t))_{P,Q}$ is upper-triangular and $c_{PP}(t) = 1$ by
the above mentioned result.

Let $(c^{QP}(t))$ be the inverse matrix $(c_{QP}(t))^{-1}$. Let
\[
  u_{RP}(t) \defeq \sum_Q c^{RQ}(t^{-1}) c_{QP}(t).
\]
Let $\setbox5=\hbox{A}
\overline{\rule{0mm}{\ht5}\hspace*{\wd5}}$ be the involution on
$\Z[t,t^{-1}]$ given by $t^{\pm 1} \mapsto t^{\mp 1}$.

\begin{Lemma}[\protect{Lusztig~\cite[7.10]{Lu:can}}]
There exists a unique solution $Z_{QP}(t)\in\Z[t^{-1}]$ \rom($Q\le
P$\rom) of
\begin{gather}
   Z_{RP}(t) = \sum_{Q: R \le Q \le P} \overline{Z_{RQ}(t)} u_{QP}(t),
   \label{eq:Lus1}\\
   Z_{PP}(t) = 1, \quad Z_{QP}(t) \in t^{-1}\Z[t^{-1}]
   \;\text{for $Q < P$}. \label{eq:Lus2}
\end{gather}
\end{Lemma}

This lemma is proved by induction, and holds in a general
setting. Lusztig has been using this (or its variant) in many places.

\begin{Theorem}\label{thm:mul}
The multiplicity $[M_P : L_Q]$ of a simple module $L_Q$ in a standard
module $M_P$ is equal to
\(
   Z_{QP}(1)
\).
\end{Theorem}

The proof will be given in \subsecref{subsec:prf}.

\begin{Example}\label{exm:A1''}
Let $\g = A_1$, $M_a = \operatorname{ev}_a^*(M)$ where
$M$ is the $2$-dimensional simple $\Ue$-module as before.
By steps~1,2 explained below\footnote{or direct calculation for the
definition~\eqref{eq:tana}}
we have
\begin{equation}\label{eq:exp}
   \chi_{q,t}(M_{a\ve^2}\otimes M_a) = 
   Y_{1,a}Y_{1,a\ve^2} + Y_{1,a}Y_{1,a\ve^4}^{-1}
   + Y_{1,a\ve^2}^{-1}Y_{1,a\ve^4}^{-1} 
   + t^{-1} 1.
\end{equation}
Let $P(u) = (1-au)(1-a\ve^2 u)$ (i.e., $M_P = M_{a\ve^2}\otimes M_a$),
$Q(u) = 1$ (i.e. $M_Q = $ trivial module).
Then the above algorithm gives us
\(
   Z_{QP}(t) = t^{-1}.
\)
\end{Example}

\section{Step 1}

\subsection{Some definitions}

Let $M_P$ be a standard module. Let $m_P$ be the monomial
corresponding to the {\it l\/}-highest weight vector. Let $M_P(Q/R)$
be an {\it l\/}-weight space as before. We denote by $m_{Q/R}$ the
corresponding monomial. We define
\(
  w_{k,a}(P), v_{k,a}(Q/R,P)\in\Z_{\ge 0}, u_{k,a}(Q/R)\in\Z
\)
by
\begin{gather*}
   m_P = \prod_{k\in I, a\in\C^*} Y_{k,a}^{w_{k,a}(P)},
\\
   m_{Q/R} = m_P \prod_{k\in I, a\in\C^*} A_{k,a}^{-v_{k,a}(Q/R,P)}
   = \prod_{k\in I, a\in\C^*} Y_{k,a}^{u_{k,a}(Q/R)}.
\end{gather*}

Suppose two standard modules $M_{P^1}$, $M_{P^2}$ and {\it l\/}-weight
spaces $M_{P^1}(Q^1/R^1)\subset M_{P^1}$, $M_{P^2}(Q^2/R^2)\subset
M_{P^2}$ are given. We define
\begin{equation}
\begin{split}
   & d(Q^1/R^1,P^1; {Q^2/R^2},{P^2}) \\
   \defeq \; & \sum_{k,a} \left(
   v_{k,a}(Q^1/R^1,P^1) u_{k,a\ve^{-1}}(Q^2/R^2)
   + w_{k,a\ve}(P^1) v_{k,a}(Q^2/R^2,P^2)\right).
\end{split}
\end{equation}

We also define
\begin{equation}\label{eq:dim}
   d(Q/R,P) \defeq d(Q/R,P;Q/R,P)
%   \sum_{k,a}
%   v_{k,a}(Q/R,P)\left( u_{k,a\ve^{-1}}(Q/R) + w_{k,a\ve}(P)\right)
.
\end{equation}
We denote $d(Q^1/R^1,P^1; {Q^2/R^2},{P^2})$ also by
$d(m_{Q^1/R^1},m_{P^1}; m_{Q^2/R^2},m_{P^2})$.

We need the following modification of $\chi_{q,t}$. Write
$\chi_{q,t}(M_P) = \sum_m a_m(t)\; m$, where $m$ is a monomial and
$a_m(t)$ is its coefficient. Let
\begin{equation}\label{eq:modify}
   \widetilde{\chi_{q,t}}(M_P) \defeq
   \sum_m t^{d(m,m_P)} a_m(t)\; m,
\end{equation}
where $d(m,m_P)$ is defined in \eqref{eq:dim}.\footnote{In fact,
$d(m,m_P)$ is determined from $a_m(t)$ so that $t^{d(m,m_P)} a_m(t)$
is a polynomial in $t$ with nonzero constant term.}

\subsection{}
Frenkel-Mukhin \cite[5.1,5.2]{FM} proved that the image of the
$q$-character $\chi_q$ is contained in
\begin{equation*}
  \bigcap_{k\in I}
  \left( \Z[Y_{l,a}^\pm]_{l\neq k, a\in\C^*}
    \otimes\Z[Y_{k,b}(1+A_{k,b\ve}^{-1})]_{b\in\C^*}\right).
\end{equation*}

We have the $t$-analogue of this result, replacing
$(1+A_{k,b\ve}^{-1})^n$ by
\begin{equation*}
   %\boxed{1+A_{k,b\ve}^{-1}}^{n}_t
   \left(1+A_{k,b\ve}^{-1}\right)^{n}_t
   \defeq \sum_{r=0}^n t^{r(n-r)} \begin{bmatrix} n \\ r
  \end{bmatrix}_t A_{k,b\ve}^{-r},
\end{equation*}
where $\left[\begin{smallmatrix} n \\ r
  \end{smallmatrix}\right]_t$ is the $t$-binomial
coefficient.
More precisely, we have
\begin{Theorem}\label{thm:ind}
\rom{(1)} For each $k\in I$, $\widetilde{\chi_{q,t}}(M_P)$ is
expressed as a linear combination of
\begin{equation*}
   \prod_i Y_{k,b_i}^{n_i}\left(1+A_{k,b_i\ve}^{-1}\right)^{n_i}_t
   = 
   Y_{k,b_1}^{n_1}\left(1+A_{k,b_1\ve}^{-1}\right)^{n_1}_t
   Y_{k,b_2}^{n_2}\left(1+A_{k,b_2\ve}^{-1}\right)^{n_2}_t
   \cdots
\end{equation*}
with coefficients in
\(
%\begin{equation*}
   \Z[t][Y_{l,a}^\pm]_{l\neq k, a\in\C^*}
%\end{equation*}
\),
where $b_i\in\C^*$, $n_i\in \Z_{>0}$ with $b_i\neq b_j$ for
$i\neq j$.

\rom{(2)} If $L_P$ is an {\it l\/}-fundamental representation \rom(and hence
$M_P = L_P$\rom), then $\chi_{q,t}(M_P)$ contains no {\it l\/}-dominant
monomials other than $m_P$ and the condition above uniquely determines
$\chi_{q,t}(M_P)$. 
\end{Theorem}

\begin{Remark}
The statement (1) for $t=1$ was proved by
Frenkel-Mukhin~\cite{FM}. And the proof of (2) is the same for $t=1$
and the general case, as illustrated in the following examples. In
this sense, (2) should also be creditted to them.
\end{Remark}

\subsection{Graph}

We give few examples of $\chi_{q,t}$ of {\it l\/}-fundamental
representations determined by the above theorem.

We attach to each standard module $M_P$, an oriented colored graph
$\Gamma_P$. (It is a slight modification of the graph in
\cite[5.3]{FR}.) The vertices are monomials in
$\chi_{q,t}(M_P)$. We draw an colored edge $\xrightarrow{k,a}$ from
$m_1$ to $m_2$ if $m_2 = m_1 A_{k,a}^{-1}$. We also write the
multiplicity of the monomials in $\chi_{q,t}(M_P)$.

\begin{Example}\label{exm:A3}
Let $\g = A_3 = \algsl_4$ and $M_P = L(\Lambda_2)_1$. Then the corresponding 
graph $\Gamma_P$ is
\begin{equation*}
\begin{CD}
   Y_{2,1} @>{2,\ve}>> Y_{1,\ve}Y_{2,\ve^2}^{-1} Y_{3,\ve}
   @>{{1,\ve^2}}>> Y_{1,\ve^3}^{-1} Y_{3,\ve} @. @. @. \\
   @. @V{{3,\ve^2}}VV @VV{{3,\ve^2}}V @. @. \\
   @. Y_{1,\ve} Y_{3,\ve^3}^{-1} @>{{1,\ve^2}}>>
   Y_{1,\ve^3}^{-1} Y_{2,\ve^2} Y_{3,\ve^3}^{-1} @>{{2,\ve^3}}>>
   Y_{2,\ve^4}^{-1}.
\end{CD}
\end{equation*}
Let us explain how we determine this graph inductively. We start with
the {\it l\/}-highest weight $Y_{2,1}$. We know that its coefficient
is $1$. Applying \thmref{thm:ind}(1) with $k=2$, we get
$Y_{1,\ve}Y_{2,\ve^2}^{-1} Y_{3,\ve}$ with coefficient $1$. Then we
apply \thmref{thm:ind}(1) with $k=1$ to get $Y_{1,\ve^3}^{-1}
Y_{3,\ve}$. And so on. All multiplicities are $1$ in this case.
\end{Example}

For $\g = A_n$, it is known that the coefficients of
$\chi_{q,t}(L(\Lambda_k)_a)$ are all $1$.\footnote{
More generally, if the coefficients of $\alpha_k$ in the highest root
is $1$, then the same holds. This result easily follows from the
theory of quiver varieties. {\bf Exercise}: Check this using the above 
algorithm.
}
Thus $\chi_{q,t}(L(\Lambda_k)_a) = \chi_{q,t=1}(L(\Lambda_k)_a)$.

\begin{Example}\label{exm:D4}
Let $\g = D_4$ and $M_P = L(\Lambda_2)_{1}$. The graph $\Gamma_P$ is
Figure~\ref{fig:D4}.
\begin{figure}[p]
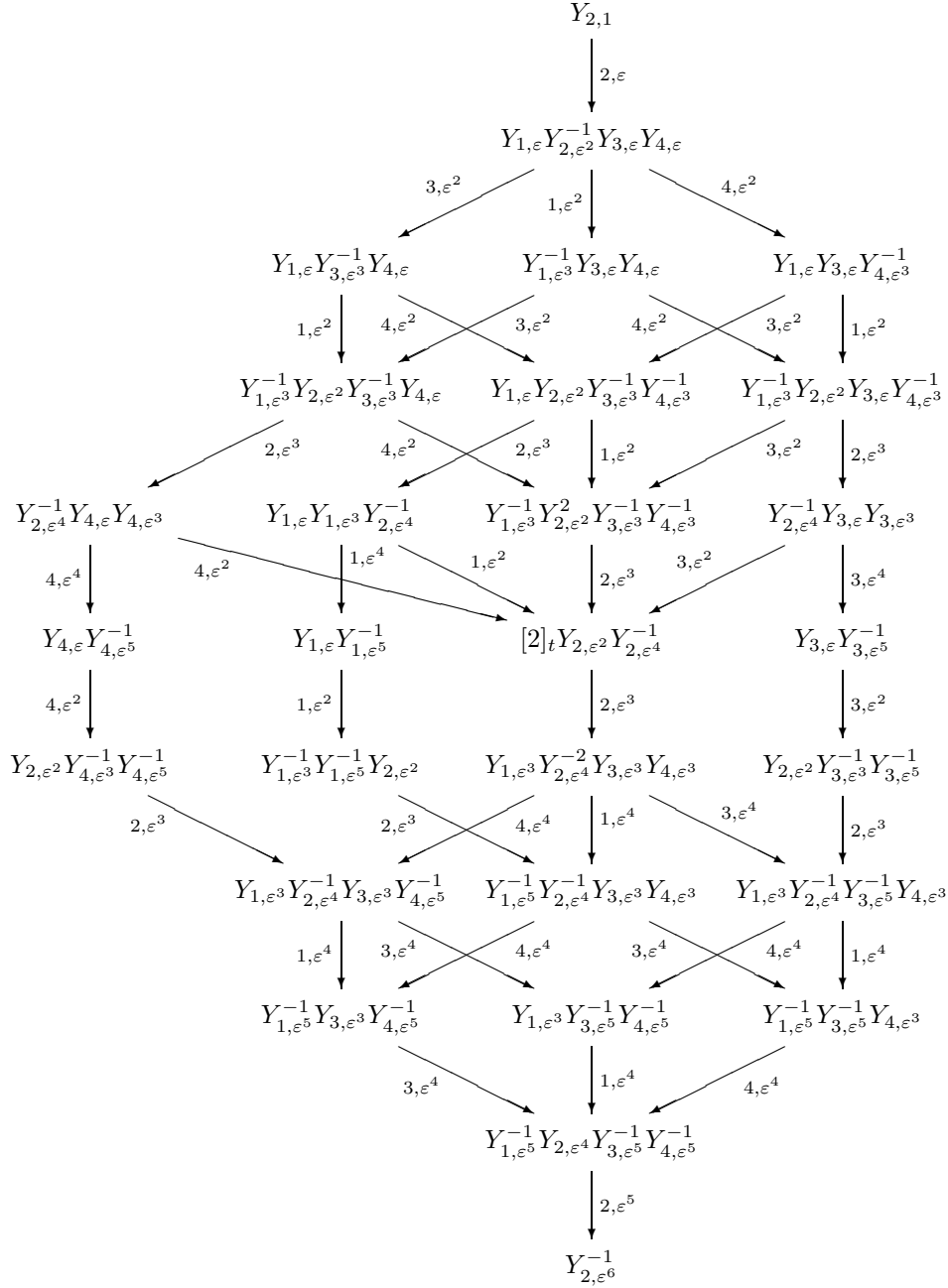

\begin{center}
\leavevmode
\begin{equation*}
\divide\dgARROWLENGTH by 3
\dgARROWPARTS=6
\begin{diagram}
  \node[3]{Y_{2,1}} \arrow{s,r}{2,\ve}
\\
  \node[3]{Y_{1,\ve}Y_{2,\ve^2}^{-1}Y_{3,\ve}Y_{4,\ve}}
  \arrow{sw,t,3}{3,\ve^2} \arrow{s,l,3}{1,\ve^2}
  \arrow{se,t,3}{4,\ve^2}
\\
  \node[2]{Y_{1,\ve} Y_{3,\ve^3}^{-1} Y_{4,\ve}}
  \arrow{s,l}{1,\ve^2} \arrow{se,b,1}{4,\ve^2}
  \node{Y_{1,\ve^3}^{-1} Y_{3,\ve} Y_{4,\ve}}
  \arrow{sw,b,1}{3,\ve^2}\arrow{se,b,1}{4,\ve^2}
  \node{Y_{1,\ve} Y_{3,\ve} Y_{4,\ve^3}^{-1}}
  \arrow{sw,b,1}{3,\ve^2} \arrow{s,r}{1,\ve^2}
\\
  \node[2]{Y_{1,\ve^3}^{-1} Y_{2,\ve^2} Y_{3,\ve^3}^{-1} Y_{4,\ve}}
  \arrow{sw,b,1}{2,\ve^3} \arrow{se,b,1}{4,\ve^2}
  \node{Y_{1,\ve} Y_{2,\ve^2} Y_{3,\ve^3}^{-1} Y_{4,\ve^3}^{-1}}
  \arrow{sw,b,1}{2,\ve^3}
  \arrow{s,r}{1,\ve^2}
  \node{Y_{1,\ve^3}^{-1} Y_{2,\ve^2} Y_{3,\ve} Y_{4,\ve^3}^{-1}}
  \arrow{sw,b,1}{3,\ve^2}
  \arrow{s,r}{2,\ve^3}
\\
  \node{Y_{2,\ve^4}^{-1}Y_{4,\ve}Y_{4,\ve^3}}
  \arrow{s,l}{4,\ve^4}
  \arrow{ese,b,1}{4,\ve^2}
  \node{Y_{1,\ve} Y_{1,\ve^3} Y_{2,\ve^4}^{-1}}
  \arrow{s,r,1}{1,\ve^4}
  \arrow{se,t,3}{1,\ve^2}
  \node{Y_{1,\ve^3}^{-1} Y_{2,\ve^2}^2 Y_{3,\ve^3}^{-1}
  Y_{4,\ve^3}^{-1}}
  \arrow{s,r}{2,\ve^3}
  \node{Y_{2,\ve^4}^{-1} Y_{3,\ve} Y_{3,\ve^3}}
  \arrow{sw,t,3}{3,\ve^2}
  \arrow{s,r}{3,\ve^4}
\\
  \node{Y_{4,\ve}Y_{4,\ve^5}^{-1}}
  \arrow{s,l}{4,\ve^2}
  \node{Y_{1,\ve}Y_{1,\ve^5}^{-1}}
  \arrow{s,l}{1,\ve^2}
  \node{[2]_t Y_{2,\ve^2}Y_{2,\ve^4}^{-1}}
  \arrow{s,r}{2,\ve^3}
  \node{Y_{3,\ve}Y_{3,\ve^5}^{-1}}
  \arrow{s,r}{3,\ve^2}
\\
  \node{Y_{2,\ve^2}Y_{4,\ve^3}^{-1}Y_{4,\ve^5}^{-1}}
  \arrow{se,b,1}{2,\ve^3}
  \node{Y_{1,\ve^3}^{-1}Y_{1,\ve^5}^{-1} Y_{2,\ve^2}}
  \arrow{se,b,1}{2,\ve^3}
  \node{Y_{1,\ve^3}Y_{2,\ve^4}^{-2}Y_{3,\ve^3}Y_{4,\ve^3}}
  \arrow{sw,b,1}{4,\ve^4}\arrow{s,r,2}{1,\ve^4}
  \arrow{se,t,3}{3,\ve^4}
  \node{Y_{2,\ve^2}Y_{3,\ve^3}^{-1}Y_{3,\ve^5}^{-1}}
  \arrow{s,r}{2,\ve^3}
\\
  \node[2]{Y_{1,\ve^3}  Y_{2,\ve^4}^{-1} Y_{3,\ve^3} Y_{4,\ve^5}^{-1}}
  \arrow{s,l}{1,\ve^4}
  \arrow{se,b,1}{3,\ve^4}
  \node{Y_{1,\ve^5}^{-1} Y_{2,\ve^4}^{-1} Y_{3,\ve^3} Y_{4,\ve^3}}
  \arrow{sw,b,1}{4,\ve^4}\arrow{se,b,1}{3,\ve^4}
  \node{Y_{1,\ve^3}  Y_{2,\ve^4}^{-1}Y_{3,\ve^5}^{-1}Y_{4,\ve^3}}
  \arrow{sw,b,1}{4,\ve^4} \arrow{s,r}{1,\ve^4}
\\
  \node[2]{Y_{1,\ve^5}^{-1} Y_{3,\ve^3} Y_{4,\ve^5}^{-1}}
  \arrow{se,b,2}{3,\ve^4}
  \node{Y_{1,\ve^3} Y_{3,\ve^5}^{-1} Y_{4,\ve^5}^{-1}}
  \arrow{s,r}{1,\ve^4}
  \node{Y_{1,\ve^5}^{-1} Y_{3,\ve^5}^{-1} Y_{4,\ve^3}}
  \arrow{sw,b,2}{4,\ve^4}
\\
  \node[3]{Y_{1,\ve^5}^{-1}Y_{2,\ve^4}Y_{3,\ve^5}^{-1}Y_{4,\ve^5}^{-1}}
  \arrow{s,r}{2,\ve^5}
\\
  \node[3]{Y_{2,\ve^6}^{-1}}
\end{diagram}
\end{equation*}
\caption{The graph for $L(\Lambda_2)_{1}$}
\label{fig:D4}
\end{center}
\end{figure}
It is known that the restriction of $M_P$ to a $\Ue$-module is a
direct sum of the adjoint representation and the trivial
representation. This fact is reflected in $\chi_{q,t}(M_P)$ where
$Y_{2,\ve^2} Y_{2,\ve^4}^{-1}$ has the coefficient $[2]_t$ and all
others has $1$. Note that the number of monomials is $28$, which is
the dimension of the adjoint representation. See also
Example~\ref{exm:D4'} below.
\end{Example}

Let us give a more complicated example.
\begin{Example}
Let $\g = A_2$ and $M_P 
= L(\Lambda_2)_{\ve}^{\otimes 2}\otimes L(\Lambda_1)_1$.
Although this is {\it not\/} an {\it l\/}-fundamental representation,
$\chi_{q,t}(M_P)$ has no {\it l\/}-dominant terms other than $m_P$,
so the condition
\thmref{thm:ind}(1) gives us $\chi_{q,t}$. The graph is
Figure~\ref{fig:A2}.

\begin{figure}[htbp]
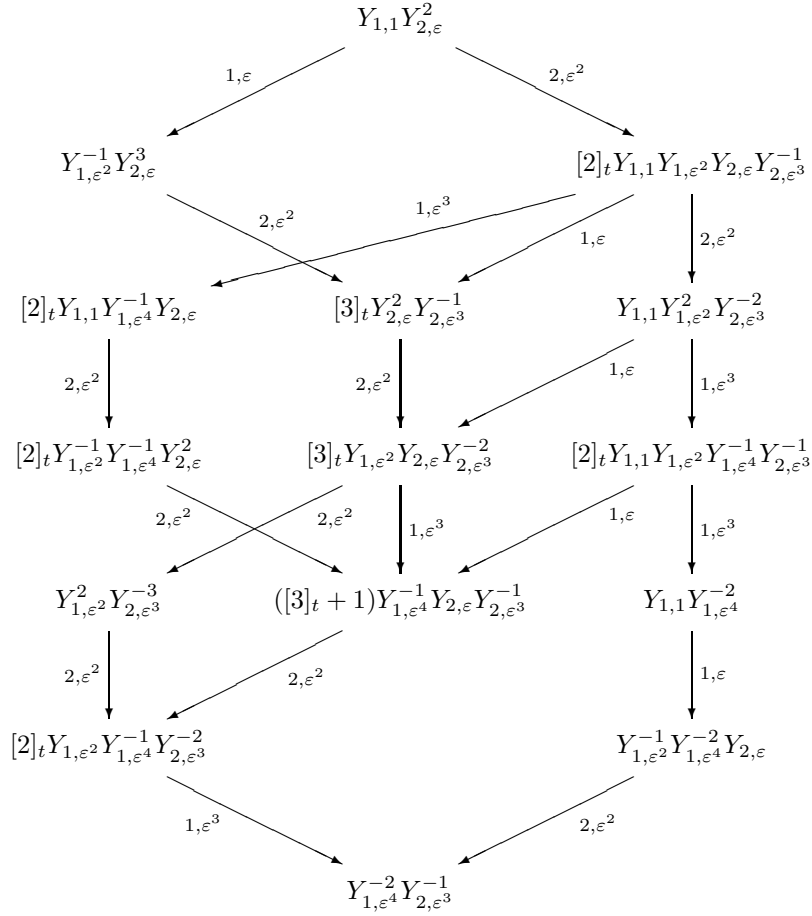

\begin{center}
\leavevmode
\begin{equation*}
\divide\dgARROWLENGTH by 2
\dgARROWPARTS=6
\begin{diagram}
\node[2]{Y_{1,1} Y_{2,\ve}^2}
\arrow{sw,t,3}{1,\ve} \arrow{se,t,3}{2,\ve^2}
\\
\node{Y_{1,\ve^2}^{-1} Y_{2,\ve}^3}
\arrow{se,t,3}{2,\ve^2}
\node[2]{[2]_t Y_{1,1} Y_{1,\ve^2} Y_{2,\ve} Y_{2,\ve^3}^{-1}}
\arrow{wsw,t,2}{1,\ve^3}
\arrow{sw,b,2}{1,\ve}
\arrow{s,r,3}{2,\ve^2}
\\
\node{[2]_t Y_{1,1} Y_{1,\ve^4}^{-1} Y_{2,\ve}}
\arrow{s,l,3}{2,\ve^2}
\node{[3]_t Y_{2,\ve}^2 Y_{2,\ve^3}^{-1}}
\arrow{s,l,3}{2,\ve^2}
\node{Y_{1,1} Y_{1,\ve^2}^2 Y_{2,\ve^3}^{-2}}
\arrow{sw,b,1}{1,\ve}
\arrow{s,r,3}{1,\ve^3}
\\
\node{[2]_t Y_{1,\ve^2}^{-1} Y_{1,\ve^4}^{-1} Y_{2,\ve}^2}
\arrow{se,b,1}{2,\ve^2}
\node{[3]_t Y_{1,\ve^2} Y_{2,\ve} Y_{2,\ve^3}^{-2}}
\arrow{sw,b,1}{2,\ve^2}
\arrow{s,r,3}{1,\ve^3}
\node{[2]_t Y_{1,1} Y_{1,\ve^2} Y_{1,\ve^4}^{-1} Y_{2,\ve^3}^{-1}}
\arrow{sw,b,1}{1,\ve}
\arrow{s,r,3}{1,\ve^3}
\\
\node{Y_{1,\ve^2}^2 Y_{2,\ve^3}^{-3}}
\arrow{s,l,3}{2,\ve^2}
\node{([3]_t + 1)Y_{1,\ve^4}^{-1} Y_{2,\ve} Y_{2,\ve^3}^{-1}}
\arrow{sw,b,2}{2,\ve^2}
\node{Y_{1,1} Y_{1,\ve^4}^{-2}}
\arrow{s,r,3}{1,\ve}
\\
\node{[2]_t Y_{1,\ve^2} Y_{1,\ve^4}^{-1} Y_{2,\ve^3}^{-2}}
\arrow{se,b,2}{1,\ve^3}
\node[2]{Y_{1,\ve^2}^{-1} Y_{1,\ve^4}^{-2} Y_{2,\ve}}
\arrow{sw,b,2}{2,\ve^2}
\\
\node[2]{Y_{1,\ve^4}^{-2} Y_{2,\ve^3}^{-1}}
\end{diagram}
\end{equation*}
\caption{The graph for $L(\Lambda_2)_{\ve}^{\otimes 2}\otimes L(\Lambda_1)_1$}
\label{fig:A2}
\end{center}
\end{figure}
\end{Example}

\begin{Remark}
As we can see in above examples, the crystal graphs are subgraphs of
$\Gamma_P$. The set of vertices is the same, but the set of arrows is
smaller. We would like to discuss this further elsewhere.
\end{Remark}

\section{Step 2}

\subsection{}

Let $M_P = M(s^1; (k_{\alpha_1}^1), (n_{\alpha_1}^1)) \otimes
M(s^2; (k_{\alpha_2}^2), (n_{\alpha_2}^2)) \otimes \cdots$ be a standard
module with $s^i/s^j\notin \ve^{\Z}$ as in \thmref{thm:VV}.
\begin{Proposition}
We have
%\(
\begin{equation*}
   \chi_{q,t}(M_P) = \chi_{q,t}
   (M(s^1; (k_{\alpha_1}^1)_{\alpha_1}, (n_{\alpha_1}^1)_{\alpha_1}))
   \chi_{q,t}
   (M(s^2; (k_{\alpha_2}^2)_{\alpha_2}, (n_{\alpha_2}^2)_{\alpha_2})) \cdots
\end{equation*}
%\)
if $s^i /s^j\notin \ve^\Z$ for $i\neq j$.
\end{Proposition}

Thus it is enough to study
\begin{equation*}
  \chi_{q,t}(M(s; (k_\alpha)_\alpha, (n_\alpha)_\alpha))
  = \chi_{q,t}(L(\Lambda_{k_1})_{\ve^{n_1}s}\otimes
   L(\Lambda_{k_2})_{\ve^{n_2}s}\otimes\cdots).
\end{equation*}

Let
\begin{equation*}
   \chi_{q,t}(L(\Lambda_{k_\alpha})_{\ve^{n_\alpha}s}))
   = \sum_{r_\alpha} a_{m_\alpha,r_\alpha}(t) m_{\alpha,r_\alpha},
\end{equation*}
where $m_{\alpha,r_\alpha}$ is a monomial in
$Y_{k,a}^\pm$ and $a_{m_\alpha,r_\alpha}(t)\in\Z[t,t^{-1}]$ is its
coefficient.

If $t=1$, $\chi_{q,1}$ is a ring homomorphism, hence we have
\begin{equation*}
  \chi_{q,1}(M(a; (k_\alpha)_\alpha, (n_\alpha)_\alpha))
  = \sum_{r_1,r_2,\dots} \prod_\alpha
  a_{m_r,r_\alpha}(1) m_{\alpha,r_\alpha}.
\end{equation*}

\begin{Theorem}
Let $P^\alpha$ be the Drinfeld polynomial of
$L(\Lambda_{k_\alpha})_{\ve^{n_\alpha}s}$.
Then we have
\begin{equation*}
  \chi_{q,t}(M(a; (k_\alpha), (n_\alpha)))
  = \sum_{r_1,r_2,\dots} 
  t^{
  \sum_{\alpha, \beta} \pm d(m_{\alpha,r_\alpha}, m_{P^\alpha};
  m_{\beta,r_\beta},m_{P^\beta})}
  \prod_\alpha
  a_{m_r,r_\alpha}(t) m_{\alpha,r_\alpha},
\end{equation*}
where the sign for $d(m_{\alpha,r_\alpha}, m_{P^\alpha};
m_{\beta,r_\beta},m_{P^\beta})$ is $-$ if $\alpha \le \beta$ and $+$
otherwise.
\end{Theorem}

\begin{Example}\label{exm:A1}
For $\g = A_1$, we have
\[
  d(Y_{1,a\ve^2}^{-1}, Y_{1,a}; Y_{1,a}, Y_{1,a}) = 1,
\qquad d(Y_{1,a}, Y_{1,a}; Y_{1,a}^{-1}, Y_{1,a\ve^{-2}}) = 1
\]
and all others are $0$.
Then we get \eqref{eq:exp}.
%\begin{equation*}
%   \chi_{q,t}(M_{a\ve^2}\otimes M_a) = Y_{1,a} Y_{1,a\ve^2}
%   + Y_{1,a} Y_{1,a\ve^4}^{-1} + Y_{1,a\ve^2}^{-1} Y_{1,a\ve^4}^{-1}
%   + t^{-1} Y_{1,a\ve^2} Y_{1,a\ve^2}^{-1}.
%\end{equation*}

If $P = (1-au)^n$, we get
\begin{equation*}
  \chi_{q,t}(M_P) = \sum_{r=0}^n
  \begin{bmatrix} n \\ r\end{bmatrix}_t 
  Y_{1,a}^{n-r} Y_{1,a\ve^2}^{-r}
\end{equation*}
from $\chi_{q,t}(L(\Lambda_1)_a) = Y_{1,a} + Y_{1,a\ve^2}^{-1}$.
This also follows directly from the definition~\eqref{eq:tana}
below. The $t$-binomial coefficients appear as Poincar\'e polynomials
of Grassmann manifolds.
\end{Example}

\section{Restrition to $\Ue$}\label{sec:rest}

Finite dimensional simple $\Ue$-modules are classified by highest
weights.
Let $\operatorname{Res}M_P$ be the restriction of a standard module
$M_P$ to a $\Ue$-module. It decomposes into a sum of various
simple modules.
Once $\chi_{q}(M_P)$ is computed, the character of
$\operatorname{Res}M_P$ is given by replacing $Y^\pm_{k,a}$ by
$y_k^\pm$ (\thmref{thm:FR}(2)). Combining with the knowledge of
characters of simple finite dimensional $\Ue$-modules, we can
determine the multiplicity of simple modules in
$\operatorname{Res}M_P$.

Characters of simple finite dimensional $\Ue$-modules are the same as
that of simple $\g$-modules, hence are known. However, we express them
in terms of $\chi_{q,t}$ in this section.

\subsection{}
For a dominant weight $\bw = \sum w_k\Lambda_k$ we denote by $L_\bw$
the simple highest weight $\Ue$-module with the highest weight $\bw$.

We consider a standard module $M_P$ with $\deg P_k = w_k$.
By the `continuity' of $M_P$ on $P$, $\operatorname{Res}M_P$ depends
only on $w_k = \deg P_k$, and not on $P$ itself. Let us denote the
multiplicity of $L_{\bw'}$ in $\operatorname{Res}M_P$ by
$Z_{\bw',\bw}$, i.e.,
\begin{equation*}
   \operatorname{Res}M_P
    = \bigoplus_{\bw'} L_{\bw'}^{\oplus Z_{\bw',\bw}}.
\end{equation*}

We will give a formula expressing $Z_{\bw',\bw}$ in terms of
$\chi_{q,t}(M_P)$. Although we can give algorithm for arbitary
$P$ in principle, the following choice will make the formula simple.

Choose and fix orientations of edges in the Dynkin diagram. We define
integer $m(k)$ for each vertex $k$ so that $m(k) - m(l) = 1$ if we
have an oriented edge from $k$ to $l$, i.e., $k\to l$. Then we define
$P$ by
\begin{equation*}
   P_k(u) = (1 - u\ve^{m(k)})^{w_k}.
\end{equation*}

Let $\widetilde{\chi_{q,t}}(M_P)$ as in \eqref{eq:modify}.
Let $\widetilde{\chi_t}(M_P)\in\Z[t]\otimes \Z[t_k^\pm]_{k\in I}$ be
a $t$-analogue of the ordinary character which is obtained from
$\widetilde{\chi_{q,t}}(M_P)$ by sending $Y_{k,a}^\pm$ to $y_k^\pm$.

For another dominant weight $\bw' = \sum w'_k\Lambda_k$, let
\begin{equation*}
   c_{\bw',\bw}(t) \defeq \text{the coefficient of $\prod y_k^{w'_k}$ in
   $\widetilde{\chi_{t}}(M_P)$}.
\end{equation*}
The matrix $(c_{\bw',\bw}(t))_{\bw',\bw}$ is upper-triangular with
respect to the usual order on weights, and diagonal entries are all
$1$.

\begin{Theorem}\label{thm:wmul}
$c_{\bw',\bw}(0)$ is the weight multiplicity of $\bw'$ in the highest
weight module $L_\bw$ with the highest weight $\bw$.
\end{Theorem}

This is just a simple rephrasing of a main result in \cite{Na:1994,Na:1998}.
The proof will be given in \subsecref{subsec:wmul}.

Note that $c_{\bw',\bw}(1)$ gives the weight multiplicity of $\bw'$ in
$\operatorname{Res}M_P$ since $\widetilde{\chi_{t=1}}$ is the ordinary
character. Thus we have
\begin{equation*}
   c_{\bw'',\bw}(1) = \sum_{\bw'} c_{\bw'',\bw'}(0)\; Z_{\bw',\bw}.
\end{equation*}
This equation determines the multiplicity $Z_{\bw',\bw}$ only from the 
knowledge of $\chi_{q,t}$.

According to a conjecture of Lusztig~\cite{Lu:ferm} together with a
formula \eqref{eq:Poin} below, $c_{\bw',\bw}(t)$ should be written by
ferminonic form of Hatayama el al.~\cite{HKOTY}. More precisely, we
should have
%\begin{equation*}
\(
   \sum_{\bw'} c^{\bw'',\bw'}(0)\, c_{\bw',\bw}(t) = M(\bw,\bw'',t^2),
\)
%\end{equation*}
where $(c^{\bw'',\bw'}(0))$ is the inverse matrix of
$(c_{\bw'',\bw'}(0))$. See \cite{Lu:ferm} for the defintion of
$M(\bw,\bw',q)$.
Although this formula can be checked in many examples, the complexity
of the combinatorics prevent us from proving it in full generality.
Conjecturally $M(\bw,\bw',q=1)$ gives us the multiplicities of the
restriction of $M_P$ (Kirillov-Reshetikhin\footnote{In fact, they
consider more general modules, not necessarily standard modules.}).
Thus the conjecture is compatible with our result in this section.

\subsection{}
\begin{Example}\label{exm:A1'}
Let $\g = A_1$ and $\bw = 2\Lambda_1$. We take $P = (1-u)^2$ by the
above choice.
By Example~\ref{exm:A1}, we
have
\begin{equation*}
   \widetilde{\chi_{t}}(M_P) = y_1^2 + (1+t^2) + y_1^{-2}.
\end{equation*}
Thus $Z_{0,\bw} = 1$. Since $\operatorname{Res}(M_P) =
L_{\Lambda_1}\otimes L_{\Lambda_1} = L_{2\Lambda_1}\oplus L_0$, this is the
correct answer~!
\end{Example}

\begin{Example}
Let $\g = A_3$ and $\bw = \Lambda_2$. By Example~\ref{exm:A3} all the
coefficients of $\chi_{q,t}(L(\Lambda_2)_1)$ are $1$. Hence
$\operatorname{Res}M_P = \operatorname{Res} L(\Lambda_2)_1$ is simple
as a $\Ue$-module.
\end{Example}

\begin{Example}\label{exm:D4'}
Let $\g = D_4$, $\bw = \Lambda_2$, $\bw' = 0$. By Example~\ref{exm:D4} 
we have $c_{\bw',\bw}(t) = 4 + t^2$. Thus $Z_{\bw',\bw} = 1$, i.e.
$\operatorname{Res}(L(\Lambda_2)_1) = L_{\Lambda_2}\oplus L_0$.
\end{Example}

\section{Quiver varieties}\label{sec:quiver}

In this section, we give the definition of $\chi_{q,t}$ and prove
Theorems~\ref{thm:mul},~\ref{thm:wmul}.
As we mentioned, those proofs are essentially given in
\cite{Na:1994,Na:1998} and \cite{Na-qaff} respectively. The only
things we do here are translation of results into the language of
$\chi_{q,t}$.
We believe that this section gives good introductions to
\cite{Na:1994,Na:1998,Na-qaff}.

\subsection{}\label{subsec:classical}

Let $\bw = \sum w_k \Lambda_k$ ($w_k\in\Z_{\ge 0}$) be a dominant
weight of the finite dimensional Lie algebra $\g$.
In \cite{Na:1994,Na:1998,Na-qaff}, we have attached to each $\bw$, a
map $\pi\colon \Mw\to \M_0(\infty,\bw)$ with the following properties:
\begin{enumerate}
\item $\Mw$ is a finite disjoint union of nonsingular quasi-projective
  varieties of various dimensions.
\item $\M_0(\infty,\bw)$ is an affine algebraic variety.
\item $\pi$ is a projective morphism.
\item There exist actions of $G_\bw\times\C^*$ on $\Mw$ and
$\M_0(\infty,\bw)$ such that $\pi$ is equivariant.
\item $\M_0(\infty,\bw)$ is a cone, and the vertex (denoted by $0$) is
the unique fixed point of the $\C^*$-action (restriction of
$G_\bw\times\C^*$-action to the second factor).
\end{enumerate}
Here $G_\bw = \prod_{k\in I} \GL(w_k,\C)$.

We consider the fiber product
\begin{equation*}
   \Zw \defeq \M(\bw)\times_{\M_0(\infty,\bw)} \M(\bw).
\end{equation*}
The convolution product makes the (Borel-Moore) homology group
$H_*(\Zw,\C)$ into an associative (noncommutative) algebra.
One of main results in \cite{Na:1998} is a construction of a
surjective algebra homomorphism
\begin{equation*}
   \bU(\g) \to H_{\topdeg}(\Zw,\C),
\end{equation*}
where $\bU(\g)$ is the universal enveloping algebra of $\g$ ({\bf NB}:
not a `quantum' version). Here $H_{\topdeg}(\ )$ means the degree
$=\dim_\R\Zw$ part of the homology group. More precisely, we take
degree = dimension part on each connected components of $\Zw$, and
then make the direct sum. Note that the the dimension differs on
various components.

Let $\Law = \pi^{-1}(0)$. It is known that $\Mw$ has a holomorphic
symplectic form such that $\Law$ is a lagrangian subvariety. The
convolution makes $H_{\topdeg}(\Law,\C)$ (the top degree part of the
Borel-Moore homology group, in the same sense as above) into an
$H_{\topdeg}(\Zw,\C)$-module. It is a $\bU(\g)$-module by the above
homomorphism. By \cite[10.2]{Na:1998} it is the simple finite dimensional
$\bU(\g)$-module $L_\bw$ with highest weight $\bw$. And connected
components $\M(\bv,\bw)$ of $\Mw$ are parametrized by vectors $\bv =
\sum v_k \alpha_k$ ($\alpha_k$ is the $k$th simple root of $\g$) so that
\begin{equation*}
   H_{\topdeg}(\Law,\C) = \bigoplus_\bv
   H_{\topdeg}(\M(\bv,\bw)\cap\Law,\C)
\end{equation*}
is the weight space decomposition of the simple highest weight module
$L_\bw$, where $H_{\topdeg}(\M(\bv,\bw)\cap\Law,\C)$ has weight $\bw -
\bv$. In particular, $\bv=0$ corresponds to the highest weight vector.
In fact, $\M(0,\bw)$ is consisting of a single point.

The space $\M_0(\infty,\bw)$ has a stratification
\begin{equation*}
   \M_0(\infty,\bw) = \bigcup \Mreg(\bv,\bw),
\end{equation*}
where $\bv$ runs over the set of vectors such that $\bw-\bv$ is a
weight of $L_\bw$ which is {\it dominant\/} \cite[\S3]{Na:1994}.

\subsection{}

Let us give the $\Ul$-version of the construction of the previous
subsection.

We use the following notation: Let $R(G)$ denote the representation
ring of a linear algebraic group $G$. If $G$ acts a quasi-projective
variety $X$, $K^G(X)$ denotes the Grothendieck group of
$G$-equivariant coherent sheaves on $X$.

The representation ring $R(G_\bw\times\C^*)$ of $G_\bw\times\C^*$ is
isomorphic to the tensor product $R(G_\bw)\otimes_\Z
R(\C^*)$. Moreover, $R(\C^*)$ is isomorphic to $\Z[q,q^{-1}]$, where
$q$ is the canonical $1$-dimensional representation of $\C^*$.

The convolution makes the Grothendieck group
$K^{G_\bw\times\C^*}(\Zw)$ into a $R(G_\bw\times\C^*) =
R(G_\bw)[q,q^{-1}]$-algebra. One of main results in \cite{Na-qaff} is
a construction of an algebra homomorphism
\begin{equation*}
   \Uli\otimes_{\Z} R(G_\bw)\to
   K^{G_\bw\times\C^*}(\Zw)/\operatorname{torsion}.
\end{equation*}

By the equivariance of $\pi$, $\Law = \pi^{-1}(0)$ is invariant under
$G_\bw\times\C^*$. The convolution makes $K^{G_\bw\times\C^*}(\Law)$
into a $K^{G_\bw\times\C^*}(\Zw)$-module. Moreover, it is free of
finite rank over $R(G_\bw\times\C^*)$ \cite[\S7]{Na-qaff}. It is a
$\Uli\otimes_\Z R(G_\bw)$-module by the above homomorphism. By
\cite[\S13]{Na-qaff}, it contains a vector $m_0$ such that
\begin{equation}
\label{eq:unistd}
\begin{gathered}
  e_{k,r}\ast m_0 = 0, \qquad
  \left(\Uli^-\otimes_\Z R(G_\bw)\right) \ast m_0
      = K^{G_\bw\times\C^*}(\Law),
\\
  \psi^{\pm}_{k}(z)\ast m_0 
  = q^{w_k}
  \left(\frac{\Wedge_{-1/qz} q^{-1} W_k}
  {\Wedge_{-q/z} q^{-1} W_k}\right)^\pm \ast m_0
    \quad\text{for $k\in I$}.
\end{gathered}
\end{equation}
The right hand side of the third equation needs an explanation: First
$W_k$ is the vector representation of $\GL(w_k,\C)$, considered as a
$G_\bw\times\C^*$-module. Then $\Wedge_u V = \sum u^i \Wedge^i V$.
Since $\Wedge_{-q/z} q^{-1} W_k$ is $1 - (1/z) W_k + \dots$ ($1$ is
the trivial module), we can define $\left(\Wedge_{-q/z} q^{-1}
  W_k\right)^{-1}$ as a formal power series in $1/z$. This gives us
the case $(\ )^+$ of the above formula. In the case $(\ )^-$, we
expand as
\(
  \Wedge_{-q/z} q^{-1} W_k = 
  (-1/z)^{w_k}
  \left(\Wedge^{w_k} W_k - z \Wedge^{w_k-1} W_k +
    \cdots\right)
\).
Then $\Wedge^{w_k} W_k$ is an invertible element, we can also define
$\left(\Wedge_{-q/z} q^{-1} W_k\right)^{-1}$.
The vector $m_0$ is the canonical generator of
$K^{G_\bw\times\C^*}(\M(0,\bw))$. (Recall $\M(0,\bw)$ is a point.)

The module $K^{G_\bw\times\C^*}(\Law)$ should be considered as a
`universal' standard module since standard modules are obtained from
it by specializations as we explain now.

Let $a = (s,\ve)\in G_\bw\times\C^*$ be a semisimple element. It
defines a homomorphism $\chi_a\colon R(G_\bw\times\C^*)\to \C$ by
sending a representation to the value of the character at $a$. Then
\begin{equation}\label{eq:std}
   K^{G_\bw\times\C^*}(\Law) \otimes_{R(G_\bw\times\C^*)} \C
\end{equation}
is a module over $\Ule = \Uli\otimes_{\Z[q,q^{-1}]}\C$.
By \eqref{eq:unistd} it is a finite-dimensional {\it l\/}-highest
weight module. This is the {\it standard module\/} $M_P$, where
$P_k(u)  = \chi_a(\Wedge_{-u}\linebreak[0] q^{-1} W_k)$.
Note that the set of conjugacy classes of $a = (s,\ve)$ bijectively
corresponds to the set of $I$-tuple of polynomials $P$ with $\deg P_k = w_k$.

\subsection{}\label{subsec:qchar}

Let $A$ be the Zariski closure of $a^\Z$ in $G_\bw\times\C^*$. It is
an abelian reductive group. We have
$K^{G_\bw\times\C^*}(\Law)\otimes_{R(G_\bw\times\C^*)} R(A) \cong
K^A(\Law)$ \cite[\S7]{Na-qaff}. Since $\chi_a$ factors through $R(A)$,
the standard module $M_P$ is isomorphic to
$K^A(\Law)\otimes_{R(A)}\C$. By Thomason's localization theorem, it is
isomorphic to $K(\Law^A)\otimes_{\Z}\C$, where $\Law^A$ is the fixed
point set. Furthermore, it is isomorphic to $H_*(\Law^A,\C)$ via the
Chern character homomorphism \cite[\S7]{Na-qaff}.

Let $\pi^A\colon\Mw^A\to\M_0(\infty,\bw)^A$ be the restriction of the
map $\pi\colon\Mw\to\M_0(\infty,\bw)$ to the fixed point set. Let
$\Mw^A = \bigsqcup_\rho \M(\rho)$ be the decomposition into connected
components. Each $\M(\rho)$ is a nonsingular quasi-projective variety.
Then we have the direct sum decomposition
\begin{equation*}
   M_P \cong H_*(\Law^A,\C) \cong \bigoplus_\rho H_*(\M(\rho)\cap\Law,\C).
\end{equation*}
In \cite[\S13, \S14]{Na-qaff} we have shown that this is the {\it
l\/}-weight space decomposition of $M_P$. In particular, the index
$\rho$ can be considered as an {\it l\/}-weight of $M_P$. Thus we have
arrived at a geometric interpretation of $\chi_q$:
\begin{equation*}
   \chi_q(M_P) 
   = \sum_\rho \dim H_*(\M(\rho)\cap\Law,\C)\; m_\rho,
\end{equation*}
where $m_\rho$ is the monomial corresponding to the {\it l\/}-weight $\rho$.

Now we define the $t$-analogue $\chi_{q,t}$ by
\begin{equation}\label{eq:tana}
   \chi_{q,t}(M_P)
   \defeq
   \sum_\rho \sum_k \dim H_k(\M(\rho)\cap\Law,\C)\;
   t^{k-\dim_\C \M(\rho)} \, m_\rho.
\end{equation}

By \cite[\S14]{Na-qaff} we have a stratification
\begin{equation*}
   \M_0(\infty,\bw)^A = \bigcup_\rho \Mreg(\rho),
\end{equation*}
consisting of nonsingular locally closed subvarieties.
Here the index set $\{ \rho\}$ is the subset of the above index set
consisting of {\it l-dominant\/} {\it l\/}-weights.

\subsection{\protect{Proof of \thmref{thm:mul}}}\label{subsec:prf}

The {\it l\/}-highest weight $P$ is fixed throughout the proof. Thus
the dominant weight vector $\bw$ and the element $a =
(s,\ve)\in G_\bw$ are fixed.

We change the notation now. If $\rho$ corresponds to an {\it
l\/}-weight space $M_P(Q/R)$, we denote above $\M(\rho)$ by
$\M(Q/R,P)$. We also denote by $\Mreg(Q,P)$ for above $\Mreg(\rho)$ if
$\rho$ corresponds to an {\it l\/}-dominant {\it l\/}-weight $Q$. Thus
we have
\begin{equation*}
   \M(\bw)^A = \bigsqcup_{Q/R} \M(Q/R,P), \quad
   \M_0(\infty,\bw)^A = \bigcup_Q \Mreg(Q,P).   
\end{equation*}
In this notation $H_*(\M(P,P)\cap\Law,\C)$ is the {\it l\/}-highest
weight space. Since $\M(0,\bw)$ is a single point as we explained, we
have $\M(P,P) = \M(0,\bw)$. We also have $\Mreg(P,P) = \{ 0\}$.

\begin{Lemma}\label{lem:slice}
\rom{(1)} $\dim_\C \M(Q/R,P) = d(Q/R,P)$. $\dim_\C \Mreg(Q,P) = d(Q,P)$.

\rom{(2)} If $\Mreg(Q,P)\subset \overline{\Mreg(R,P)}$, then
$R\le Q$.

\rom{(3)} Choose $x\in\Mreg(Q,P)$. Then 
$(\pi^A)^{-1}(x)\cap\M(S/T,P)$ is isomorphic to
$\M(S/T,Q)\cap\Law$.
\end{Lemma}

\begin{proof}
(1) The first equation is the dimension formula
\cite[4.1.6]{Na-qaff}. The second equation follows from $\dim_\C \Mreg(Q,P) =
\dim_\C \M(Q,P)$, which is clear from the definition \cite[\S4]{Na-qaff}.

(2),(3) The results are known or trivial for $Q = P$. Now use the
transversal slice at $x\in\Mreg(Q,P)$ \cite[\S3]{Na-qaff} to reduce
a general case to this case.
\end{proof}

Let $D^b(\M_0(\infty,\bw)^A)$ be the bounded derived category of
complexes of sheaves such that cohomology sheaves are constant along
each stratum $\Mreg(Q,P)$. Let $IC(\Mreg(Q,P))$ be the intersection
homology complex associated with the constant local system
$\C_{\Mreg(Q,P)}$ on $\Mreg(Q,P)$. By using the transversal slice
\cite[\S3]{Na-qaff}, one can check that it is an object in
$D^b(\M_0(\infty,\bw)^A)$. Let $\C_{\M(Q/R,P)}$ be the constant local
system on $\M(Q/R,P)$. Then $\pi^A_*(\C_{\M(Q/R,P)})$ is an object of
$D^b(\M_0(\infty,\bw)^A)$ again by the transversal slice
argument. Using the decomposition theorem of
Beilinson-Bernstein-Deligne, we have shown that there exists an
isomorphism in $D^b(\M_0(\infty,\bw)^A)$:
\begin{equation}\label{eq:decomp}
   \pi^A_*(\C_{\M(R,P)}[\dim_\C\M(R,P)]) 
        \cong \bigoplus_{Q,k} L_{Q,k}(R,P)\otimes IC(\Mreg(Q,P))[k]
\end{equation}
for some vector space $L_{Q,k}(R,P)$ \cite[14.3.2]{Na-qaff}.
Since $\pi^A(\M(R,P))\subset\overline{\Mreg(R,P)}$ by
definition~\cite[\S4]{Na-qaff}, the summation runs over $Q\ge R$ by
\lemref{lem:slice}.
Let
\begin{equation*}
   L_{RQ}(t) \defeq
   \sum_{k} \dim L_{Q,k}(R,P)\, t^{-k}.
\end{equation*}
Applying the Verdier duality to the both hand side of
\eqref{eq:decomp} and using the self-duality of 
$\pi^A_*(\C_{\M(R,P)}[\dim_\C\M(R,P)])$ and $IC(\Mreg(Q,P))$, we find
$L_{RQ}(t) = \overline{L_{RQ}(t)}$.

Choose a point $x_Q$ from $\Mreg(Q,P)$ for each stratum. Let
$i_{x_Q}\colon \{ x_Q\} \to \M_0(\infty,\bw)^A$ denote the
inclusion. Consider
\begin{equation*}
   H^{k}(i_{x_Q}^! \pi^A_* \C_{\M(R,P)}[\dim \M(R,P)])
   = H_{\dim_\C \M(R,P)- k}((\pi^A)^{-1}(x_Q)\cap \M(R,P),\C).
\end{equation*}
By \lemref{lem:slice}(3) this is isomorphic to
\(
%\begin{equation*}
   H_{\dim_\C \M(R,P)- k}(\M(R,Q)\cap\Law,\C)
%\end{equation*}
\). Therefore we have
\begin{equation}\label{eq:C_RQ}
\begin{split}
   & \sum_k \dim H^{k}(i_{x_Q}^! \pi^A_* \C_{\M(R,P)}[\dim \M(R,P)])\,
   t^{\dim_\C \M(Q,P)-k}
\\
  = \; &
   \sum_d \dim H_{d}(\M(R,Q)\cap\Law,\C)\,
   t^{d + \dim_\C \M(Q,P) - \dim_\C \M(R,P)}
  = c_{RQ}(t),
\end{split}
\end{equation}
where we used $\dim_\C\M(R,P) - \dim_\C\M(Q,P) = \dim_\C\M(R,Q)$ in the
last equality.

By \cite[14.3.10]{Na-qaff}, we have
\begin{equation*}
   [M_Q:L_R] = \dim H^*(i_{x_Q}^! IC(\Mreg(R,P))).
\end{equation*}
(In fact, we defined the standard module $M_Q$ as
$H_*((\pi^A)^{-1}(x_Q),\C)$ in \cite[\S13]{Na-qaff}, which apriori
depends on $P$. Thus the definition coincides only when
$Q=P$. However, by using the transversal slice, we can show that the
right hand side is the same for both definitions. cf.\
\lemref{lem:slice}.)

Let
\begin{equation*}
   Z_{RQ}(t) \defeq
   \sum_k \dim H^k(i_{x_Q}^! IC(\Mreg(R,P)))\, t^{\dim\Mreg(Q,P)-k}.
\end{equation*}
We have $[M_Q:L_R] = Z_{RQ}(1)$. By the defining property of the
intersection homology, $Z_{RQ}(t)$ satisfies \eqref{eq:Lus2}.

Substituting \eqref{eq:decomp} into \eqref{eq:C_RQ}, we get
\begin{equation*}
   c_{SQ}(t) = \sum_{R} L_{SR}(t) Z_{RQ}(t).
\end{equation*}
Now $L_{SR}(t) = \overline{L_{SR}(t)}$ implies \eqref{eq:Lus1}.
This completes the proof of \thmref{thm:mul}.

\subsection{Proof of \protect{\thmref{thm:wmul}}}\label{subsec:wmul}

By the result explained in \subsecref{subsec:classical}, the weight
multiplicity of $\bw'$ in $L_\bw$ is equal to
\begin{equation*}
   \dim H_{\topdeg}(\M(\bw - \bw',\bw)\cap \Law,\C).
\end{equation*}
The assertion follows from more general formula
\begin{equation}\label{eq:Poin}
   \widetilde{\chi_{t}}(M_P)
   = \sum_{\bw'} 
   \sum_d H_d(\M(\bw-\bw',\bw)\cap \Law,\C)\,
      t^{\dim_\C \M(\bw-\bw',\bw) - d}\;
   \prod_k y_k^{w'_k}.
\end{equation}
Note that $\M(\bw-\bw',\bw)\cap \Law$ is a lagrangian subvariety in
$\M(\bw-\bw',\bw)$, so we have $\topdeg = \dim_\C\M(\bw-\bw',\bw)$.

In order to prove \eqref{eq:Poin}, we use \cite[5.7]{Na:1994}, where
the Betti numbers are given in terms of those of fixed point
components. It looks almost the same as above. However, there is one
significant difference. The $\C^*$-action used there is different from
our $\C^*$-action used here, defined in \cite[\S2]{Na-qaff}.
This is the reason why we choose $P$ and corresponding $a = (s,\ve)$
as explained in \secref{sec:rest}. Then $A = \overline{a^{\Z}}$ is
isomorphic to $\C^*$ and the action is the same as the $\C^*$-action
considered in \cite[\S5]{Na:1994}.

We decompose $\Mw^A = \bigsqcup \M(\rho)$ into connected components as
before. By \cite[5.7]{Na:1994} we have\footnote{
In fact, this formula even holds for general $P$ if we
replace $\dim_\C \M(\bw-\bw',\bw)-d$ by a suitable degree. However,
this degree shift is given by a complicated expression in $\rho$. So
our choice of $P$ is most economical.
}
\begin{equation*}
   \dim H_d(\M(\bw-\bw',\bw)\cap \Law,\C)
   = \sum_\rho \dim H_{\dim_\C \M(\bw-\bw',\bw)-d}(\M(\rho),\C),
\end{equation*}
where the summation runs over the set of
$\rho$ such that the corresponding monomial $m_\rho$ is sent to
$\prod_k y_k^{w'_k}$ after $Y_{k,a} \to y_k$.
The $\C^*$-action makes $\M_0(\infty,\bw)^{\C^*} = \{ 0\}$, so
$\M(\rho) = \M(\rho)\cap\Law$.
Hence the above expression coincides with the definition
of the coefficients $\widetilde{\chi_{t}}$.

\subsection*{Acknowledgement}
We would like to thank E.~Frenkel and E.~Mukhin for explanations of
their algorithm computing $q$-characters of fundamental
representations.

%\bibliographystyle{amsplain}
%\bibliography{qchar}

\begin{thebibliography}{10}

\bibitem{CP-book}
V.~Chari and A.~Pressley, 
\emph{A guide to quantum groups}, Cambridge University Press,
  Cambridge, 1994.

\bibitem{CP-rep}
\bysame, \emph{Quantum affine algebras and their representations},
  Representations of groups (Banff, AB, 1994), Amer. Math. Soc., Providence,
  RI, 1995, pp.~59--78.

\bibitem{CP-roots}
\bysame, \emph{Quantum affine algebras at roots of unity}, Represent. Theory
  \textbf{1} (1997), 280--328 (electronic).

\bibitem{CP:Weyl}
\bysame, \emph{{Integrable and {W}eyl modules for
  quantum affine $sl_2$}}, preprint, arXiv:math.QA/0007123.

\bibitem{FM}
E.~Frenkel and E.~Mukhin, \emph{{Combinatorics of q-characters of
  finite-dimensional representations of quantum affine algebras}},
preprint, arXiv:math.QA/9911112.

\bibitem{FR}
E.~Frenkel and N.~Reshetikhin, \emph{The $q$-characters of
representations of quantum affine algebras and deformations of
$\mathscr {W}$-algebras}, Recent developments in quantum affine
algebras and related topics (Raleigh, NC, 1998), Contemp. Math., 248,
Amer. Math. Soc., Providence, RI, 1999, pp.~163--205.

\bibitem{HKOTY}
G.~Hatayama, A.~Kuniba, M.~Okado, T.~Takagi and Y.~Yamada,
\emph{Remarks on fermionic formula}, Recent developments in quantum
affine algebras and related topics (Raleigh, NC, 1998),
Contemp. Math., 248, Amer. Math. Soc., Providence, RI, 1999,
pp.~243--291.

\bibitem{Kas}
M.~Kashiwara, \emph{Crystal bases of modified quantized enveloping
  algebra}, Duke Math. J. \textbf{73} (1994), no.~2, 383--413.

\bibitem{KL}
D.~Kazhdan and G.~Lusztig, \emph{Representations of {C}oxeter groups and
  {H}ecke algebras}, Invent. Math. \textbf{53} (1979), no.~2, 165--184.

\bibitem{Lu:can}
G.~Lusztig, \emph{Canonical bases arising from quantized enveloping algebras},
  J. Amer. Math. Soc. \textbf{3} (1990), no.~2, 447--498.

\bibitem{Lu:ferm}
\bysame, \emph{{{F}erminonic form and {B}etti numbers}}, preprint,
arXiv:math.QA/0005010.

\bibitem{Na:1994}
H.~Nakajima, \emph{Instantons on {A}{L}{E} spaces, quiver varieties, and
  {K}ac-{M}oody algebras}, Duke Math. J. \textbf{76} (1994), no.~2, 365--416.

\bibitem{Na-hom}
\bysame, \emph{Homology of moduli spaces instantons on {A}{L}{E} spaces~{I}},
J. Differential Geom. \textbf{40} (1994), 105--127.

\bibitem{Na:1998}
\bysame, \emph{Quiver varieties and {K}ac-{M}oody algebras}, Duke Math. J.
  \textbf{91} (1998), no.~3, 515--560.

\bibitem{Na-qaff}
\bysame, \emph{{Quiver varieties and finite dimensional representations
  of quantum affine algebras}}, preprint, arXiv:math.QA/9912158,
to appear in J. Amer. Math. Soc.

\bibitem{VV-std}
M.~Varagnolo and E.~Vasserot, \emph{{Standard modules of quantum affine
  algebras}}, preprint, arXiv:math.QA/0006084.

\end{thebibliography}

\providecommand{\bysame}{\leavevmode\hbox to3em{\hrulefill}\thinspace}

\end{document}